\documentstyle[12pt]{article}
\def\Bbb R{{\rm \bf R}}
\def\proclaim#1{\vskip2mm{\bf #1}\em}
\def\endproclaim{\em \vskip2mm}
\def\tag#1{\eqno(#1)}
\def\gathered{\begin{array}{c}}
\def\endgathered{\end{array}}
\def\text{\mbox}
\begin{document}

\title {Inverse obstacle scattering problems with a single incident wave and
the logarithmic differential of the indicator function
in the Enclosure Method}
\author{Masaru IKEHATA\footnote{
Department of Mathematics,
Graduate School of Engineering,
Gunma University, Kiryu 376-8515, JAPAN,
ikehata@math.sci.gunma-u.ac.jp}
}
\date{April 3, 2011}
\maketitle
\begin{abstract}
This paper gives a remark on the Enclosure Method by considering
inverse obstacle scattering problems with a single incident wave
whose governing equation is given by the Helmholtz equation in two dimensions.
It is concerned with the indicator function in the
Enclosure Method.  The previous indicator function is essentially
real-valued since only its absolute value is used.
In this paper, another method for the use of the indicator function is introduced.
The method employs the {\it logarithmic differential} with respect to the independent
variable of the indicator function and yields directly the coordinates
of the vertices of the convex hull of unknown polygonal sound-hard obstacles
or thin ones.  The convergence rate of the formulae is better than that
of the previous indicator function.  Some other applications of this method
are also given.

\noindent
AMS: 35R30


\end{abstract}


\section{Introduction and statements of the main results}

The aim of this paper is to add further new knowledge on the {\it Enclosure Method}.
The Enclosure Method was originally introduced in
\cite{IE, IE3} for inverse boundary value problems for elliptic
equations.  The method aims at extracting information about the
location and shape of unknown discontinuity embedded in a known
reference medium that gives an effect on the propagation of the
signal, such as an obstacle, inclusion, crack, etc. 
The method can be divided into two versions.
One is a version that employs infinitely many pairs of input and output data,
that is, the Dirichlet-to-Neumann map (or Neumann-to-Dirichlet map).
Another is a version that employs a single set of input and output data.
We call this second version the {\it single measurement version} of the Enclosure Method.
We refer the reader to \cite{IW, IP, IH} 
for recent applications of the single measurement version of 
the Enclosure Method.

This paper is concerned with the {\it indicator function} in the
single measurement version of the Enclosure Method. 
In \cite{IE4}, having the single measurement version of the Enclosure Method,
the author considered the reconstruction issue of inverse
obstacle scattering problems of acoustic wave in two dimensions.
The problem is to reconstruct a two dimensional obstacle from the
Cauchy data on a circle surrounding the obstacle of the total wave
field generated by a {\it single} incident {\it plane wave} with a
{\it fixed} wave number $k>0$. The author established
an extraction formula of the value of the support function
at a generic direction which yields information about the {\it
convex hull} of {\it polygonal} sound-hard obstacles. 
However, the indicator function used in \cite{IE4} is essentially real-valued since 
only its absolute value is used.
In this paper, another method for the use
of the indicator function is introduced.
It is shown that the {\it
logarithmic differential} with respect to the independent variable
of the indicator function yields directly the coordinates of the
vertices of the convex hull of unknown polygonal sound-hard
obstacles or thin ones.

Let us describe our main results.
First consider a {\it polygonal} obstacle denoted by $D$, that is:
$D\subset\Bbb R^2$ takes the form $D_1\cup\cdots\cup D_m$ with
$1\le m<\infty$ where each $D_j$ is open and a polygon; $\overline
D_j\cap\overline D_{j'}=\emptyset$ if $j\not=j'$.

The total wave field $u$ outside obstacle $D$ takes the form $u(x;d,k)=e^{ikx\cdot d}+w(x)$
with $k>0$, $d\in S^1$
and satisfies
$$\begin{array}{c}
\displaystyle
\triangle u+k^2 u=0\,\,\text{in}\,\Bbb R^2\setminus\overline D,\\
\\
\displaystyle
\frac{\partial u}{\partial\nu}=0\,\,\text{on}\,\partial D,\\
\\
\displaystyle
\lim_{r\longrightarrow\infty}\sqrt{r}\left(\frac{\partial w}{\partial r}-ikw\right)=0,
\end{array}
$$
where $r=\vert x\vert$ and $\nu$ denotes the unit outward normal relative to $D$.
The last condition above is called the Sommerfeld radiation condition.
Some further information about $u$
are in order.  $u$ belongs to $C^{\infty}(\Bbb R^2\setminus\overline D)$ and
satisfies $u\vert_{B}\in H^1((\Bbb R^2\setminus\overline D)
\cap B)$ for a large open disc $B$ containing $\overline D$.
This restricts the singularity of $u$ at the corner of $D$.
The boundary condition for $\partial u/\partial\nu$ on $\partial D$
means that $D$ is a {\it sound-hard} obstacle and
should be considered as a weak sense.

Note that, all the results in this paper $k$ is just positive and there is {\it no other assumption} on $k$.

Let $B_R$ be an open disc with radius $R$ centered at a fixed point
satisfying $\overline D\subset B_R$.
We assume that $B_R$ is {\it known}.  Our data are  $u=u(\,\cdot\,;d,k)$ and $\partial u/\partial\nu$ on $\partial B_R$
for a fixed $d$ and $k$,
where $\nu$ is the unit outward normal relative to $B_R$.
Let $\omega$ and $\omega^{\perp}$ be two unit
vectors perpendicular to each other.  We always choose the orientation of $\omega^{\perp}$ and $\omega$ coincides with
$\mbox{\boldmath $e$}_1$ and $\mbox{\boldmath $e$}_2$ and thus $\omega^{\perp}$ is unique.

Recall the {\it support function} of $D$:  $h_D(\omega)=\sup_{x\in D}x\cdot\omega$.
We say that $\omega$ is {\it regular} with respect to $D$ if the set
$\displaystyle
\partial D\cap\{x\in\Bbb R^2\,\vert\,x\cdot\omega=h_D(\omega)\}$
consists of only one point.

Set $c_{\tau}(\omega)=\tau\,\omega+i\sqrt{\tau^2+k^2}\,\omega^{\perp}$ with $\tau>0$.
Let $v_{\tau}(x)=e^{x\cdot c_{\tau}(\omega)}$.  This $v$ satisfies the Helmholtz equation
in the whole plane.

Define
$$\displaystyle
I(\tau;\omega, d, k)
=\int_{\partial B_R}\left(\frac{\partial u}{\partial\nu}v_{\tau}
-\frac{\partial v_{\tau}}{\partial\nu}u\right)
dS.
\tag {1.1}
$$
This complex valued-function of $\tau$ is called the {\it indicator function} in
the single measurement version of the Enclosure Method.
In \cite{IE4} Ikehata has established the formula:
$$\displaystyle
\lim_{\tau\longrightarrow\infty}\frac{1}{\tau}
\log
\left\vert I(\tau;\omega, d, k)
\right\vert
=h_D(\omega),
\tag {1.2}
$$
provided $\omega$ is regular with respect to $D$.
In this formula one makes use of only the {\it absolute
value} of indicator function $I(\tau;\omega,d,k)$.
Thus one needs two regular directions
$\omega$ for determining a single vertex of the convex hull of $D$
since formula (1.2) gives only a single line on which the vertex
lies.  Here we present a method for the use of the complex values
of the indicator function which directly yields the coordinates of a
vertex of the convex hull of $D$ with indicator functions for a
single regular direction $\omega$

Since $(\omega^{\perp})^{\perp}=-\omega$,
we have $\sqrt{\tau^2+k^2}\,\omega+i\tau\omega^{\perp}=ic_{\tau}(\omega^{\perp})$.
This gives
$$\displaystyle
\partial_{\tau}v_{\tau}
=\frac{i}{\sqrt{\tau^2+k^2}}\,x\cdot c_{\tau}(\omega^{\perp})v_{\tau}
$$
and thus we have
$$\displaystyle
\frac{\partial}{\partial\nu}(\partial_{\tau}v_{\tau})
=\frac{i}{\sqrt{\tau^2+k^2}}
\left\{c_{\tau}(\omega^{\perp})\cdot\nu+(x\cdot c_{\tau}(\omega^{\perp}))(c_{\tau}(\omega)\cdot\nu)\right\}v_{\tau}
\tag {1.3}
$$
and
$$\displaystyle
I'(\tau;\omega,d,k)
=\int_{\partial B_R}
\left(\frac{\partial u}{\partial\nu}\cdot\partial_{\tau}v_{\tau}-
\frac{\partial}{\partial\nu}(\partial_{\tau}v_{\tau})u\right)dS.
\tag {1.4}
$$

Our first result is the following theorem.

\proclaim{\noindent Theorem 1.1.}
Let $\omega$ be regular with respect to $D$.
Let $x_0\in\partial D$ be the point with $x_0\cdot\omega=h_D(\omega)$.
There exists a $\tau_0>0$ such that, for all $\tau\ge\tau_0$
$\vert I(\tau;\omega,d,k)\vert>0$
and the formula
$$\displaystyle
\lim_{\tau\longrightarrow\infty}
\frac{\displaystyle I'(\tau;\omega,d,k)}
{I(\tau;\omega,d,k)}
=h_D(\omega)+ix_0\cdot\omega^{\perp},
\tag {1.5}
$$
is valid.

\endproclaim

Theorem 1.1 means that, as $\tau\longrightarrow\infty$ the {\it
logarithmic differential} of $I(\tau;\omega,d,k)$ with respect to
$\tau$ converges to $x_0\cdot\omega+ix_0\cdot\omega^{\perp}$. The
convergence rate of (1.5) is better than that of (1.2).  For this
see Remark 2.1 in Section 2.

{\bf\noindent Remark 1.1.}
By (1.4) $I'(\tau;\omega,d,k)$ can be considered as another indicator function, however,
it is easy to see that from (1.2) and (1.5) we have
$$\displaystyle
\lim_{\tau\longrightarrow\infty}\frac{1}{\tau}
\log
\left\vert I'(\tau;\omega, d, k)
\right\vert
=h_D(\omega).
$$
Thus from this formula we can not obtain any new information about
$D$.

Here we explain why (1.5) gives further information about $D$.
Let $\omega$ be {\it regular} with respect to $D$.
We denote by $x(\omega)=(x(\omega)_1,x(\omega)_2)$ the single point
in $\displaystyle
\partial D\cap\{x\in\Bbb R^2\,\vert\,x\cdot\omega=h_D(\omega)\}$.
Since it holds that
$$\displaystyle
x\cdot(\omega+i\omega^{\perp})=(x_1-ix_2)(\omega_1+i\omega_2),\,\,x\in\Bbb R^2,
\tag {1.6}
$$
from (1.5) we have
$$\displaystyle
x(\omega)_1=\text{Re}\,\left\{(\omega_1+i\omega_2)\lim_{\tau\longrightarrow\infty}
\overline{\left(\frac{\displaystyle I'(\tau;\omega,d,k)}
{I(\tau;\omega,d,k)}\right)}\right\}
\tag {1.7}
$$
and
$$\displaystyle
x(\omega)_2=\text{Im}\,\left\{(\omega_1+i\omega_2)\lim_{\tau\longrightarrow\infty}
\overline{\left(\frac{\displaystyle I'(\tau;\omega,d,k)}
{I(\tau;\omega,d,k)}\right)}\right\}.
\tag {1.8}
$$
The set $I(D)$ of all directions which are not regular with respect to $D$ is finite.
Let $I(D)=\{(\cos\,\theta_j,\sin\,\theta_j)\,\vert\,0\le \theta_1<\cdots<\theta_N<2\pi\}$
and $k=1,2$.
Define
$$\displaystyle
x(\theta)_k=\left\{
\begin{array}{lr}
\displaystyle
x((\cos\,\theta,\sin\,\theta))_k
, & \quad\text{if $\theta\in[0,\,2\pi[\setminus\{\theta_j\,\vert\,j=1,\cdots,N\}$,}\\
\\
\displaystyle \infty, & \quad\text{if $\theta=\theta_j$, $j=1,\cdots,N$.}
\end{array}
\right.
$$
and extend it as the $2\pi$-periodic function of $\theta\in\Bbb
R$. Since $D$ is polygonal, both $x(\theta)_1$ and $x(\theta)_2$
are {\it piece-wise constant} and one of which has discontinuity
at $\theta=\theta_j$ in the following sense: for each $j$ it holds
that $x(\theta_j+0)_1\not=x(\theta_j-0)_1$ or
$x(\theta_j+0)_2\not=x(\theta_j-0)_2$. Therefore one can expect
that computing both $x(\omega)_1$ and $x(\omega)_2$ for
sufficiently many $\omega$ via formulae (1.7) and (1.8), one can
estimate $I(D)$. This is a new information extracted from the
indicator function in the Enclosure Method.

Another implication of Theorem 1.1 is the following idea.
Given $y\in\Bbb R^2$ define
$$\displaystyle
I(\tau;y,\omega,d,k)
=e^{-y\cdot(\tau\omega+i\sqrt{\tau^2+k^2}\omega^{\perp})}I(\tau;\omega,d,k).
$$
This corresponds to substitute $x\mapsto e^{(x-y)\cdot(\tau\omega+i\sqrt{\tau^2+k^2}\omega^{\perp})}$
instead of $v_{\tau}$ into (1.1).
Since
$$\displaystyle
I'(\tau;y,\omega,d,k)
=-y\cdot\left(\omega+i\frac{\tau}{\sqrt{\tau^2+k^2}}\omega^{\perp}\right)I(\tau;y,\omega,d,k)
+e^{-y\cdot(\tau\omega+i\sqrt{\tau^2+k^2}\omega^{\perp})}I'(\tau;\omega,d,k),
$$
it follows from (1.5) that
$$\displaystyle
\lim_{\tau\longrightarrow\infty}
\frac{I'(\tau;y,\omega,d,k)}
{I(\tau;y,\omega,d,k)}
=(x_0-y)\cdot(\omega+i\omega^{\perp}).
$$
This together with (1.6) yields that
$$\displaystyle
\lim_{\tau\longrightarrow\infty}
\left\vert\frac{I'(\tau;y,\omega,d,k)}
{I(\tau;y,\omega,d,k)}\right\vert
=\vert x_0-y\vert.
$$
Since $x_0$ is the unique point which minimizes the function
$y\mapsto \vert x_0-y\vert$, one possible alternative idea to find
$x_0$ is to consider the {\it minimization} problem of the
following function for a suitable $\tau$:
$$\displaystyle
y\mapsto \left\vert\frac{I'(\tau;y,\omega,d,k)}
{I(\tau;y,\omega,d,k)}\right\vert.
$$
Since this paper concentrates on only the theoretical issue of the
Enclosure Method, we leave the numerical implementation of this
idea for future research.

The result can be extended to a {\it thin}
obstacle case.
Let $\Sigma$ be the union of finitely many disjoint closed piecewise linear segments
$\Sigma_1, \Sigma_2,\cdots,\Sigma_m$.
Assume that there exists  a simply connected open set $D$
such that $D$ is a polygon and each $\Sigma_j$ consists of sides of $D$.
We assume that $\overline D\subset B_{R}$ with a $R>0$.
We denote by $\nu$ the unit outward normal on $\partial D$ relative to $B_{R}\setminus\overline D$
and set $\nu^+=\nu$ and $\nu^-=-\nu$ on $\Sigma$.
Given $k>0$ and $d\in S^1$ let $u=u(x), x\in\,\Bbb R^2\setminus\Sigma$
be the solution of the scattering problem
$$\begin{array}{c}
\displaystyle
(\triangle+k^2)u=0\,\,\text{in}\,\,\Bbb R^2\setminus\Sigma,\\
\\
\displaystyle
\frac{\partial u^{\pm}}{\partial\nu^{\pm}}=0\,\,\text{on}\,\Sigma,\\
\\
\displaystyle
\lim_{r\longrightarrow\infty}\,\sqrt{r}\left(\frac{\partial w}{\partial r}-ik\,w\right)=0,
\end{array}
$$
where $w=u-e^{ikx\cdot d}$, $u^+=u\vert_{\Bbb
R^2\setminus\overline D}$ and $u^-=u\vert_{D}$.  Note that this is
a brief description of the problem and for exact one see
\cite{IE4}. Define
$$\displaystyle
I_{\Sigma}(\tau;\omega,d,k)
=\int_{\partial B_{R}}
\left(\frac{\partial u}{\partial\nu} v_{\tau}-\frac{\partial v_{\tau}}{\partial\nu}u\right)dS.
$$

The following is our second result.

\proclaim{\noindent Theorem 1.2.}
Let $\omega$ be regular with respect to $\Sigma$.
If every end points of $\Sigma_1,\Sigma_2,\cdots,\Sigma_m$ satisfies
$x\cdot\omega<h_{\Sigma}(\omega)$, then
there exists a $\tau_0>0$ such that, for all $\tau\ge\tau_0$
$\vert I_{\Sigma}(\tau;\omega,d,k)\vert>0$
and the formula
$$\displaystyle
\lim_{\tau\longrightarrow\infty}
\frac{\displaystyle I_{\Sigma}'(\tau;\omega,d,k)}
{I_{\Sigma}(\tau;\omega,d,k)}
=h_D(\omega)+ix_0\cdot\omega^{\perp},
\tag {1.9}
$$
is valid.
\noindent
If there is an end point $x_0$ of some $\Sigma_j$ such that
$x_0\cdot\omega=h_{\Sigma}(\omega)$, then, for $d$ that is not
perpendicular to $\nu$ on $\Sigma_j$ near the point, the same
conclusions as above are valid.

\endproclaim

Note that $\nu$ on $\Sigma_j\cap B_{\eta}(x_0)$ for sufficiently small $\eta>0$
becomes a constant vector if $x_0$ is an end point of $\Sigma_j$.

In \cite{IE4} under the same assumption as Theorem 1.2, it is shown that
$$\displaystyle
\lim_{\tau\longrightarrow\infty}\frac{1}{\tau}\log\vert I_{\Sigma}(\tau;\omega,d,k)\vert
=h_{\Sigma}(\omega).
\tag {1.10}
$$
Thus (1.9) also adds a further knowledge on the use of the
indicator function in thin obstacle case.

A brief outline of this paper is as follows. Theorems 1.1 and 1.2
are proved in Sections 2 and 3, respectively. Both proofs employ
some previous computation results done in \cite{IE4} for the proof
of (1.2) and (1.10), however, some nontrivial modifications of the
computation are also required.
The idea of using the logarithmic differential of the original
indicator function developed in this paper can be applied to
several other previous applications of the Enclosure Method
published in \cite{IE3, IE100, ILAY2, ITRANS, Ik11, IP, IH}.
In Section 4 two applications of
the argument for the proof of Theorem 1.1 are given.
In Appendix
first for reader's convenience we give the proof of Proposition
2.1 which ensures an expansion of the solution of the Helmholtz
equation at a corner. The proof is focused on some technical part
that is different from the case when $k=0$. Second a proof of
Lemma 2.3 which is important for the computation of the expansion
of $I'(\tau;\omega,d,k)$ as $\tau\longrightarrow\infty$ is given.
Third a proof of some estimates that are needed for the proof of
one of two applications in Section 4 is given

Note also that in Sections 2 and 3 we simply write $v_{\tau}=v$.

\section{Proof of Theorem 1.1}

Let $x_0$ denote the single point of the set
$\{x\,\vert\,x\cdot\omega=h_D(\omega)\}\cap\partial D$. $x_0$ has
to be a vertex of $D_j$ for some $j$.
The internal angle of $D_j$ at $x_0$ is less than $\pi$ and
thus $2\pi$ minus the internal angle which we denote by $\Theta$
satisfies $\pi<\Theta<2\pi$.

In what follows we denote by
$B_R(x_0)$ the open disc with radius $R$ centered at $x_0$.
If one chooses a
sufficiently small $\eta>0$, then one can write
$$\begin{array}{c}
\displaystyle
B_{2\eta}(x_0)\cap (B_{R_1}\setminus\overline D)
=\{x_0+r(\cos\,\theta\,\text{\boldmath $a$}+\sin\,\theta\,\text{\boldmath $a$}^{\perp})\,\vert\,
0<r<2\eta,\,0<\theta<\Theta\},
\\
\\
\displaystyle
B_{\eta}(x_0)\cap\partial D=\Gamma_p\cup \Gamma_q\cup\{x_0\}
\end{array}
$$
where $\displaystyle
\text{\boldmath $a$}=\cos\,p\,\omega^{\perp}+\sin\,p\,\omega$,
$\displaystyle
\text{\boldmath $a$}^{\perp}=-\sin\,p\,\omega^{\perp}+\cos\,p\,\omega$;
$\displaystyle -\pi<p<0$;
$\displaystyle
\Gamma_p=\{x_0+r\text{\boldmath $a$}\,\vert\,0<r<\eta\}$,
$\displaystyle\Gamma_q=\{x_0+r(\cos\,\Theta\,\text{\boldmath $a$}+\sin\,\Theta\,\text{\boldmath $a$}^{\perp})\,\vert\,0<r<\eta\}$.
Note that the orientation of $\text{\boldmath $a$}$, $\text{\boldmath $a$}^{\perp}$ coincides with that of
$\text{\boldmath $e$}_1, \text{\boldmath $e$}_2$.
See also Figure 1 of \cite{IE3}.

The quantity $-p$ means the angle between two vectors $\omega^{\perp}$ and $\text{\boldmath $a$}$.
$p$ satisfies $\Theta>\pi+(-p)$. Set $q=\Theta-2\pi+p$.  Then we have $-\pi<q<p<0$ and
the expression
$$\begin{array}{c}
\displaystyle
\Gamma_p=\{x_0+r(\cos\,p\,\omega^{\perp}+\sin\,p\,\omega)\,\vert\,0<r<\eta\},\\
\\
\displaystyle
\Gamma_q=\{x_0+r(\cos\,q\,\omega^{\perp}+\sin\,q\,\omega)\,\vert\,0<r<\eta\}.
\end{array}
$$
This is the meaning of $p$ and $q$.

\noindent
We set
$$\displaystyle
u(r,\theta)=u(x),\,\,x=x_0+r(\cos\,\theta\,\text{\boldmath $a$}+\sin\,\theta\,
\text{\boldmath $a$}^{\perp}).
$$

The followng proposition describes the behaviour of $u(r,\theta)$ as $r\longrightarrow 0$.

\proclaim{\noindent Proposition 2.1(Proposition 4.2 in \cite{IE4}).}
Let $\eta$ satisfy $\eta<<1/2k$.
Then, there exists a sequence $\alpha_1,\alpha_2,\cdots,\alpha_m,\cdots$ such that:

\noindent
(1) for each $s\in]0,\,2[$
$$\displaystyle
u(r,\theta)=\sum_{m=1}^{\infty}\alpha_m J_{\mu_m}(kr)\cos \mu_m\theta,\,\,\text{in}\,
H^1(B_{s\eta}(x_0)\cap(B_R\setminus\overline D))
$$
where
$$\displaystyle
\mu_m=\frac{(m-1)\pi}{\Theta}
$$
and $J_{\mu_m}$ denotes the Bessel function of order $\mu_m$ given by the formula
$$\displaystyle
J_{\mu_m}(z)=\left(\frac{z}{2}\right)^{\mu_m}\sum_{n=0}^{\infty}
\frac{(-1)^n}{n!\Gamma(n+1+\mu_m)}\left(\frac{z}{2}\right)^{2n};
$$

\noindent
(2) as $m\longrightarrow\infty$
$$\displaystyle
\vert\alpha_m\vert=O\left(\frac{\Gamma(1+\mu_m)}{\sqrt{\mu_m}}\left(\frac{1}{\eta k}\right)^{\mu_m}\right);
$$

\noindent
(3) for each $l=1,\cdots$ there exists a positive number $C_l$ such that,
for all $r\in\,]0,\,\eta[$
$$\begin{array}{c}
\displaystyle
\vert u(r,0)-\sum_{m=1}^l\alpha_mJ_{\mu_m}(kr)\vert\le C_lr^{\mu_{l+1}}\\
\\
\displaystyle
\vert u(r,\Theta)-\sum_{m=1}^l\alpha_m(-1)^{m-1} J_{\mu_m}(kr)\vert\le C_lr^{\mu_{l+1}},\,\,0<r<\eta.
\end{array}
\tag {2.1}
$$
\endproclaim

\noindent
In \cite{IE4} the proof is omitted since it can be done along the same line
as the proof in the case when $k=0$ given in \cite{G}.  
However, there is a technical difference from the case when $k=0$
and so to make sure and for reader's convenience, 
in Appendix we give the proof which focused on the diffrence.

Let $s=\sqrt{\tau^2+k^2}+\tau$.
We have:
$$
\displaystyle
\tau=\frac{1}{2}\left(s-\frac{k^2}{s}\right),\,\,
\sqrt{\tau^2+k^2}=\frac{1}{2}\left(s+\frac{k^2}{s}\right).
$$

Recall that the proof of formula (1.2) is based on the following two facts in \cite{IE4}.

$\bullet$  As $s\longrightarrow\infty$ the {\it complete} asymptotic expansion
$$\displaystyle
e^{-i\sqrt{\tau^2+k^2}x_0\cdot\omega^{\perp}}e^{-\tau h_D(\omega)}I(\tau;\omega,d,k)
\sim
-i\sum_{n=2}^{\infty}\frac{e^{i\frac{\pi}{2}\mu_n}k^{\mu_n}\alpha_n
\{e^{ip\mu_n}+(-1)^n e^{iq\mu_n}\}}
{s^{\mu_n}},
\tag {2.2}
$$
is valid.

$\bullet$  $\exists n\ge 2$\,\, $\alpha_n\{e^{ip\mu_n}+(-1)^n e^{iq\mu_n}\}\not=0$.
Thus the quantity
$$\displaystyle
m^*=\min\{m\ge 2\,\vert\,
\alpha_m\{e^{ip\mu_m}+(-1)^me^{iq\mu_m}\}\not=0\}
$$
is well-defined.  $m^*$ depends on $k$, $d$, $D$ and $\omega$.

For the proof of Theorem 1.1 we compute the asymptotic expansion
of $I'(\tau;\omega,d,k)$ as $s\longrightarrow\infty$.
The result is:
$$\begin{array}{c}
\displaystyle
\sqrt{\tau^2+k^2}e^{-i\sqrt{\tau^2+k^2}\,x_0\cdot\omega^{\perp}}e^{-\tau h_D(\omega)}I'(\tau;\omega,d,k)\\
\\
\displaystyle
=-i\sum_{m=2}^n
\alpha_m\{e^{ip\mu_m}+(-1)^me^{iq\mu_m}\}(-\mu_m+ix_0\cdot c_{\tau}(\omega^{\perp}))
\frac{e^{i\frac{\pi}{2}\mu_m}k^{\mu_m}}{s^{\mu_m}}
+O\left(\frac{1}{s^{\mu_{n+1}}}\right).
\end{array}
\tag {2.3}
$$

The proof of (2.3) is given in Subsection 2.1.
Here we show how to prove Theorem 1.1 by assuming (2.3).

Since $\alpha_m\{e^{ip\mu_m}+(-1)^me^{iq\mu_m}\}=0$ for all $m$ with $m<m^*$ and
$\beta\equiv\alpha_{m^*}\{e^{ip\mu_{m^*}}+(-1)^{m^*}e^{iq\mu_{m^*}}\}\not=0$,
from (2.2) we have
$$\displaystyle
\lim_{\tau\longrightarrow\infty}
\tau^{\mu_{m^*}}e^{-i\sqrt{\tau^2+k^2}\,x_0\cdot\omega^{\perp}}e^{-\tau h_D(\omega)}I(\tau;\omega,d,k)
=-i\beta e^{i\frac{\pi}{2}\mu_{m^*}}\left(\frac{k}{2}\right)^{\mu_{m^*}}.
\tag {2.4}
$$
On the other hand, since
$$\displaystyle
\lim_{\tau\longrightarrow\infty}\frac{ix_0\cdot c_{\tau}(\omega^{\perp})}
{\sqrt{\tau^2+k^2}}=x_0\cdot\omega+ix_0\cdot\omega^{\perp},
$$
it follows from (2.3) that
$$
\displaystyle
\lim_{\tau\longrightarrow\infty}
\tau^{\mu_{m^*}}e^{-i\sqrt{\tau^2+k^2}\,x_0\cdot\omega^{\perp}}e^{-\tau h_D(\omega)}I'(\tau;\omega,d,k)
=-i\beta(x_0\cdot\omega+ix_0\cdot\omega^{\perp})
e^{i\frac{\pi}{2}\mu_{m^*}}\left(\frac{k}{2}\right)^{\mu_{m^*}}.
\tag {2.5}
$$
A combination of (2.4) and (2.5) ensures the validity of (1.5).
This completes the proof of Theorem 1.1.

{\bf\noindent Remark 2.1.}
Since $c_{\tau}(\omega^{\perp})$ depends on $\tau$ and thus $s$, (2.3)
is not the {\it complete} asymptotic expansion.
However, using (2.2), (2.3) and the expression
$$\displaystyle
ic_{\tau}(\omega^{\perp})
=\frac{s}{2}(\omega+i\omega^{\perp})
+\frac{k^2}{2s}(\omega-i\omega^{\perp}),
$$
one can easily obtain the following expansion:
$$\begin{array}{c}
\displaystyle
e^{-i\sqrt{\tau^2+k^2}\,x_0\cdot\omega^{\perp}}e^{-\tau h_D(\omega)}
\left\{\sqrt{\tau^2+k^2}I'
-\left(\frac{s}{2}x_0\cdot(\omega+i\omega^{\perp})
+\frac{k^2}{2s}x_0\cdot(\omega-i\omega^{\perp})\right)I\right\}
\\
\\
\displaystyle
=i\sum_{m=2}^n
\alpha_m\{e^{ip\mu_m}+(-1)^me^{iq\mu_m}\}
\frac{\mu_me^{i\frac{\pi}{2}\mu_m}k^{\mu_m}}{s^{\mu_m}}
+O\left(\frac{1}{s^{\mu_{n+1}}}\right),
\end{array}
$$
where $I=I(\tau;\omega,d,k)$ and $I'=I'(\tau;\omega,d,k)$.
Note that we have used $\pi<\Theta\le 2\pi$.

The above formula yields the second term of the expansion of the logarithmic differential
of the indicator function as $\tau\longrightarrow\infty$:
$$\displaystyle
\frac{I'(\tau;\omega,d,k)}{I(\tau;\omega,d,k)}
=h_D(\omega)+ix_0\cdot\omega^{\perp}
-\frac{\mu_{m^*}}{\tau}
+O\left(\frac{1}{\tau^2}\right).
$$
Thus the convergence rate of (1.5) is {\it better} than that of
(1.2) since (2.4) yields
$$\displaystyle
\frac{1}{\tau}\log\vert I(\tau;\omega,d,k)\vert=h_D(\omega)
-\frac{\mu_{m^*}\log\tau}{\tau}+O\left(\frac{1}{\tau}\right).
$$
From this view point one can say that (1.5) is an improvement of (1.2).
Note also that both formulae show that the accuracy of the approximation depends on
the {\it size} of $\mu_{m^*}=(m^*-1)\pi/\Theta$.

\subsection{Proof of (2.3)}

Integration by parts gives
$$\displaystyle
I'(\tau;\omega,d,k)
=\int_{\partial D} u\frac{\partial}{\partial\nu}\partial_{\tau}vdS.
$$
Localizing this integral at $x_0$, we have, modulo exponentially decaying
as $\tau\longrightarrow\infty$
$$\begin{array}{c}
\displaystyle
e^{-i\sqrt{\tau^2+k^2}\,x_0\cdot\omega^{\perp}}e^{-\tau h_D(\omega)}I'(\tau;\omega,d,k)\\
\\
\displaystyle
\sim
e^{-i\sqrt{\tau^2+k^2}\,x_0\cdot\omega^{\perp}}e^{-\tau h_D(\omega)}
\int_{\Gamma_p}
u\frac{\partial}{\partial\nu}\partial_{\tau}vdS
+e^{-i\sqrt{\tau^2+k^2}\,x_0\cdot\omega^{\perp}}e^{-\tau h_D(\omega)}\int_{\Gamma_q}
u\frac{\partial}{\partial\nu}\partial_{\tau}vdS
\\
\\
\displaystyle
\equiv
I_p(\tau)+I_q(\tau).
\end{array}
\tag {2.6}
$$

We have:$\nu=\sin\,p\,\omega^{\perp}-\cos\,p\,\omega$ and
$x=x_0+r(\cos\,p\,\omega^{\perp}+\sin\,p\,\omega)$ on $\Gamma_p$;
$\nu=-\sin\,q\,\omega^{\perp}+\cos\,q\,\omega$ and
$x=x_0+r(\cos\,q\,\omega^{\perp}+\sin\,q\,\omega)$ on $\Gamma_q$.

From those we have:

on $\Gamma_p$
$$\begin{array}{c}
\displaystyle
c_{\tau}(\omega^{\perp})\cdot\nu
=\tau\sin\,p+i\sqrt{\tau^2+k^2}\cos\,p,\\
\\
\displaystyle
c_{\tau}(\omega)\cdot\nu
=-(\tau\cos\,p-i\sqrt{\tau^2+k^2}\sin\,p),\\
\\
\displaystyle
x\cdot c_{\tau}(\omega^{\perp})
=x_0\cdot c_{\tau}(\omega^{\perp})
+r(\tau\cos\,p-i\sqrt{\tau^2+k^2}\sin\,p),\\
\\
\displaystyle
v=e^{\tau h_D(\omega)}
e^{i\sqrt{\tau^2+k^2}\,x_0\cdot\omega^{\perp}}
e^{\tau r\,\sin p}e^{i\sqrt{\tau^2+k^2}\,r\,\cos\,p};
\end{array}
\tag {2.7}
$$

on $\Gamma_q$
$$\begin{array}{c}
\displaystyle
c_{\tau}(\omega^{\perp})\cdot\nu
=-(\tau\sin\,q+i\sqrt{\tau^2+k^2}\cos\,q),\\
\\
\displaystyle
c_{\tau}(\omega)\cdot\nu
=\tau\cos\,q-i\sqrt{\tau^2+k^2}\sin\,q,\\
\\
\displaystyle
x\cdot c_{\tau}(\omega^{\perp})
=x_0\cdot c_{\tau}(\omega^{\perp})
+r(\tau\cos\,q-i\sqrt{\tau^2+k^2}\sin\,q),\\
\\
\displaystyle
v=e^{\tau h_D(\omega)}
e^{i\sqrt{\tau^2+k^2}\,x_0\cdot\omega^{\perp}}
e^{\tau r\,\sin q}e^{i\sqrt{\tau^2+k^2}\,r\,\cos\,q}.
\end{array}
\tag {2.8}
$$

It follows from (1.3) and (2.7) that
$$\begin{array}{c}
\displaystyle
\sqrt{\tau^2+k^2}I_p(\tau)\\
\\
\displaystyle
=
i\left\{(\tau\sin\,p+i\sqrt{\tau^2+k^2}\cos\,p)
-x_0\cdot c_{\tau}(\omega^{\perp})(\tau\cos\,p-i\sqrt{\tau^2+k^2}\sin\,p)\right\}\\
\\
\displaystyle
\times
\int_0^{\eta}u(r,0)
e^{\tau r\,\sin p}e^{i\sqrt{\tau^2+k^2}\,r\,\cos\,p}dr\\
\\
\\
\displaystyle
-i(\tau\cos\,p-i\sqrt{\tau^2+k^2}\sin\,p)^2
\int_0^{\eta}ru(r,0)e^{\tau r\,\sin p}e^{i\sqrt{\tau^2+k^2}\,r\,\cos\,p}dr.
\end{array}
\tag {2.9}
$$

It follows also from (1.3) and (2.8) that
$$\begin{array}{c}
\displaystyle
\sqrt{\tau^2+k^2}I_q(\tau)\\
\\
\displaystyle
=
-i\left\{(\tau\sin\,q+i\sqrt{\tau^2+k^2}\cos\,q)
-x_0\cdot c_{\tau}(\omega^{\perp})(\tau\cos\,q-i\sqrt{\tau^2+k^2}\sin\,q)\right\}\\
\\
\displaystyle
\times
\int_0^{\eta}u(r,\Theta)e^{\tau r\,\sin q}e^{i\sqrt{\tau^2+k^2}\,r\,\cos\,q}dr\\
\\
\displaystyle
+i(\tau\cos\,q-i\sqrt{\tau^2+k^2}\sin\,q)^2
\int_0^{\eta}ru(r,\Theta)e^{\tau r\,\sin q}e^{i\sqrt{\tau^2+k^2}\,r\,\cos\,q}dr.
\end{array}
\tag {2.10}
$$

Here we make use of (2.1).
Since $\alpha_1 J_{0}(k\vert x-x_0\vert)$ satisfies the Helmholtz equation
in the whole plane, the indicator function for $u$ coincides with that for $u-\alpha_1 J_0(k\vert x-\vert_0)$.
Thus one may assume, in advance $\alpha_1=0$ in the computation of the integrals in (2.9) and (2.10).

Since $p$ and $q$ satisfies $-\pi<q<p<0$, we have $\sin\,p<0$ and $\sin\,q<0$.
This gives, for $\theta=p,q$ and $\mu>0$
$$\displaystyle
\int_0^{\eta}r^{\mu}e^{\tau r\,\sin\theta}dr
=O(\tau^{-(\mu+1)}).
\tag {2.11}
$$
Set
$$\begin{array}{c}
\displaystyle
I_{\mu}(\tau,\theta)
=\int_0^{\eta}J_{\mu}(kr)e^{\tau r\,\sin\,\theta}e^{i\sqrt{\tau^2+k^2}\,r\,\cos\,\theta}dr,
\\
\\
\displaystyle
K_{\mu}(\tau,\theta)
=\int_0^{\eta}rJ_{\mu}(kr)e^{\tau r\,\sin\theta}e^{i\sqrt{\tau^2+k^2}\,r\,\cos\,\theta}dr.
\end{array}
\tag {2.12}
$$
It follows from (2.1) and (2.11) that
$$\begin{array}{c}
\displaystyle
\int_0^{\eta}u(r,0)e^{\tau r\,\sin p}e^{i\sqrt{\tau^2+k^2}\,r\,\cos\,p}dr
=\sum_{m=2}^n\alpha_mI_{\mu_m}(\tau,p)
+O(\tau^{-(\mu_{n+1}+1)}),
\end{array}
$$

$$\begin{array}{c}
\displaystyle
\int_0^{\eta}ru(r,0)e^{\tau r\,\sin p}e^{i\sqrt{\tau^2+k^2}\,r\,\cos\,p}dr
=\sum_{m=2}^n\alpha_mK_{\mu_m}(\tau,p)
+O(\tau^{-(\mu_{n+1}+2)}),
\end{array}
$$

$$\begin{array}{c}
\displaystyle
\int_0^{\eta}u(r,\Theta)e^{\tau r\,\sin q}e^{i\sqrt{\tau^2+k^2}\,r\,\cos\,q}dr
=\sum_{m=2}^n\alpha_m(-1)^{m-1}I_{\mu_m}(\tau,q)
+O(\tau^{-(\mu_{n+1}+1)}),
\end{array}
$$

$$\begin{array}{c}
\displaystyle
\int_0^{\eta}ru(r,\Theta)e^{\tau r\,\sin q}e^{i\sqrt{\tau^2+k^2}\,r\,\cos\,q}dr
=\sum_{m=2}^n\alpha_m(-1)^{m-1}K_{\mu_m}(\tau,q)
+O(\tau^{-(\mu_{n+1}+2)}).
\end{array}
$$

Substituting these into (2.9) and (2.10), we obtain

$$\begin{array}{c}
\displaystyle
\sqrt{\tau^2+k^2}I_p(\tau)\\
\\
\displaystyle
=
i\left\{(\tau\sin\,p+i\sqrt{\tau^2+k^2}\cos\,p)
-x_0\cdot c_{\tau}(\omega^{\perp})(\tau\cos\,p-i\sqrt{\tau^2+k^2}\sin\,p)\right\}\\
\\
\displaystyle
\times
\sum_{m=2}^n\alpha_mI_{\mu_m}(\tau,p)
\\
\\
\displaystyle
-i(\tau\cos\,p-i\sqrt{\tau^2+k^2}\sin\,p)^2\sum_{m=2}^n\alpha_mK_{\mu_m}(\tau,p)
+O(\tau^{-\mu_{n+1}})
\end{array}
\tag {2.13}
$$
and
$$\begin{array}{c}
\displaystyle
\sqrt{\tau^2+k^2}I_q(\tau)\\
\\
\displaystyle
=
-i\left\{(\tau\sin\,q+i\sqrt{\tau^2+k^2}\cos\,q)
-x_0\cdot c_{\tau}(\omega^{\perp})(\tau\cos\,q-i\sqrt{\tau^2+k^2}\sin\,q)\right\}\\
\\
\displaystyle
\times
\sum_{m=2}^n\alpha_m(-1)^{m-1}I_{\mu_m}(\tau,q)
\\
\\
\displaystyle
+i(\tau\cos\,q-i\sqrt{\tau^2+k^2}\sin\,q)^2
\sum_{m=2}^n\alpha_m(-1)^{m-1}K_{\mu_m}(\tau,q)
+O(\tau^{-\mu_{n+1}}).
\end{array}
\tag {2.14}
$$

In what follows we choose an arbitrary $n$ and fix.
Let $\mu=\mu_m$ with $0\le m\le n$.
We have already established that
\proclaim{\noindent Proposition 2.2(\cite{IE4}).}
As $\tau\longrightarrow\infty$
we have
$$\displaystyle
(\tau\,\cos\,\theta-i\sqrt{\tau^2+k^2}\,\sin\,\theta)I_{\mu}(\tau,\theta)
=\frac{ie^{i(\theta+\frac{\pi}{2})\mu}k^{\mu}}{s^{\mu}}
+O(s^{-\infty}).
\tag {2.15}
$$

\endproclaim

The main problem is to compute the asymptotic expansion
of the quantity
$(\tau\cos\,\theta-i\sqrt{\tau^2+k^2}\,\sin\,\theta)^2K_{\mu}(\tau,\theta)$.

We prove
\proclaim{\noindent Proposition 2.3.}
As $s\longrightarrow\infty$ we have
$$\begin{array}{c}
\displaystyle
(\tau\,\cos\,\theta-i\sqrt{\tau^2+k^2}\,\sin\,\theta)^2K_{\mu}(\tau,\theta)\\
\\
\displaystyle
=i(1-\zeta)^2
\{1+\zeta+\mu(1-\zeta)\}(1-\zeta)^{-3}
\frac{ie^{i(\theta+\frac{\pi}{2})\mu}k^{\mu}}{s^{\mu}}
+O(s^{-\infty}),
\end{array}
\tag {2.16}
$$
where $\zeta=(k/s)^2e^{2i\theta}$.

\endproclaim

For the proof we prepare the following three technical lemmas.

\noindent
The first one can be proven along the same line as Lemma 3.1 in \cite{IE4}.

\proclaim{\noindent Lemma 2.1.}
For each $l'=0,1,\cdots$ we have
$$\begin{array}{c}
\displaystyle
K_{\mu}(\tau,\theta)\\
\\
\displaystyle
=\sum_{j=0}^{l'}\frac{(-1)^j}{j!\Gamma(1+j+\mu)}
\left(\frac{k}{2}\right)^{2j+\mu}
\int_0^{\infty}r^{2j+\mu+1}e^{\tau r\sin\,\theta}
e^{i\sqrt{\tau^2+k^2}r\cos\,\theta}dr\\
\\
\displaystyle
+O\left(\frac{1}{\tau^{2+2(l'+1)+\mu}}\right).
\end{array}
\tag {2.17}
$$

\endproclaim

The second is nothing but Lemma 3.2 in \cite{IE4}.

\proclaim{\noindent Lemma 2.2.}
Let $-\pi<\theta<0$.
For each $l=0,1,\cdots$ as $s\longrightarrow\infty$
the formula
$$\begin{array}{c}
\displaystyle
\int_0^{\infty}r^{\sigma}e^{\tau r\sin\theta} e^{i\sqrt{\tau^2+k^2}r\cos\theta} dr
=\sum_{n=0}^{l}\frac{L_{\sigma,n}(\theta)}{s^{\sigma+2n+1}}
+O\left(\frac{1}{s^{\sigma+2(l+1)+1}}\right),
\end{array}
\tag {2.18}
$$
is valid where
$$\displaystyle
L_{\sigma,n}(\theta)
=ie^{i\theta}e^{i(\theta+\frac{\pi}{2})\sigma}2^{\sigma+1}(-k^2e^{2i\theta})^n\frac{\Gamma(\sigma+n+1)}
{n!}.
$$
\endproclaim

The third corresponds to Lemma 3.3 in \cite{IE4}, however,
the proof needs a careful modification of that of Lemma 3.3.
See Appendix for the proof.

\proclaim{\noindent Lemma 2.3.}
$$\displaystyle
\sum_{n_1+n_2=n}
\frac{(-1)^{n_2}\Gamma(n+2+n_2+\mu_m)}{n_1!n_2!\Gamma(1+n_2+\mu_m)}
=(-1)^n(n+1)(n+1+\mu).
\tag {2.19}
$$
\endproclaim

{\it\noindent Proof of Proposition 2.3.}
A combination of (2.17) and (2.18) for $\sigma=2j+\mu+1$ gives
$$\begin{array}{c}
\displaystyle
K_{\mu}(\tau,\theta)\\
\\
\displaystyle
=\sum_{j=0}^{l'}\frac{(-1)^j}{j!\Gamma(1+j+\mu)}
\left(\frac{k}{2}\right)^{2j+\mu}
\left\{\sum_{n=0}^{l}\frac{L_{2j+\mu+1,n}(\theta)}{s^{2j+\mu+2n+2}}
+O\left(\frac{1}{s^{2j+\mu+2(l+1)+2}}\right)\right\}\\
\\
\displaystyle
+O\left(\frac{1}{s^{2+2(l'+1)+\mu}}\right)
\\
\\
\displaystyle
=\sum_{j=0}^{l'}\sum_{n=0}^{l}\frac{(-1)^j}{j!\Gamma(1+j+\mu)}
\left(\frac{k}{2}\right)^{2j+\mu}
\frac{L_{2j+\mu+1, n}(\theta)}
{s^{2(n+j)+\mu+2}}\\
\\
\displaystyle
+O\left(\frac{1}{s^{2(l+1)+\mu+2}}\right)
+O\left(\frac{1}{s^{2(l'+1)+\mu+2}}\right).
\end{array}
\tag {2.20}
$$

Now let $l=l'$.  Then (2.20) becomes

$$\begin{array}{c}
\displaystyle
K_{\mu}(\tau,\theta)
\\
\\
\displaystyle
=
\sum_{n=0}^l
\left(\sum_{n_1+n_2=n}
\frac{(-1)^{n_2}}{n_2!\Gamma(1+n_2+\mu)}\left(\frac{k}{2}\right)^{2n_2+\mu}
L_{2n_2+\mu+1,n_1}(\theta)\right)
\times\frac{1}{s^{2n+\mu+2}}\\
\\
\displaystyle
+O\left(\frac{1}{s^{2(l+1)+\mu+2}}\right).
\end{array}
\tag {2.21}
$$

Write
$$\begin{array}{c}
\displaystyle
\sum_{n_1+n_2=n}
\frac{(-1)^{n_2}}{n_2!\Gamma(1+n_2+\mu)}\left(\frac{k}{2}\right)^{2n_2+\mu}
L_{2n_2+\mu+1,n_1}(\theta)\\
\\
\displaystyle
=\sum_{n_1+n_2=n}
\frac{(-1)^{n_2}}{n_2!\Gamma(1+n_2+\mu)}\left(\frac{k}{2}\right)^{2n_2+\mu}
\\
\\
\displaystyle
\times
ie^{i\theta}e^{i(\theta+\frac{\pi}{2})(2n_2+\mu+1)}2^{2n_2+\mu+2}(-k^2e^{2i\theta})^{n_1}
\frac{\Gamma(2n_2+\mu+n_1+2)}{n_1!}
\\
\\
\displaystyle
=2^2ie^{i\theta}e^{i(\theta+\frac{\pi}{2})(\mu+1)}k^{\mu}
(-k^2e^{2i\theta})^n
\sum_{n_1+n_2=n}
\frac{(-1)^{n_2}\Gamma(n+2+n_2+\mu)}{n_1!n_2!\Gamma(1+n_2+\mu)}.
\end{array}
\tag {2.22}
$$

From (2.19), (2.21) and (2.22) one gets
$$
\displaystyle
K_{\mu}(\tau,\theta)
=2^2ie^{i\theta}e^{i(\theta+\frac{\pi}{2})(\mu+1)}k^{\mu}
\sum_{n=0}^l\frac{(k^2e^{2i\theta})^n(n+1)(n+1+\mu)}{s^{2n+\mu+2}}
+O\left(\frac{1}{s^{2(l+1)+\mu+2}}\right).
$$
From this one gets
$$\begin{array}{c}
\displaystyle
(\tau\,\cos\,\theta-i\sqrt{\tau^2+k^2}\,\sin\,\theta)^2K_{\mu}(\tau,\theta)\\
\\
\displaystyle
=\frac{(se^{-i\theta})^2}{2^2}
(1-\zeta)^2\\
\\
\displaystyle
\times
\left\{2^2ie^{i\theta}e^{i(\theta+\frac{\pi}{2})(\mu+1)}k^{\mu}
\sum_{n=0}^l\frac{(k^2e^{2i\theta})^n(n+1)(n+1+\mu)}{s^{2n+\mu+2}}
+O\left(\frac{1}{s^{2(l+1)+\mu+2}}\right)\right\}\\
\\
\displaystyle
=e^{-2i\theta}
(1-\zeta)^2\\
\\
\displaystyle
\times
\left\{ie^{i\theta}e^{i(\theta+\frac{\pi}{2})(\mu+1)}k^{\mu}
\sum_{n=0}^l\frac{(k^2e^{2i\theta})^n(n+1)(n+1+\mu)}{s^{2n+\mu}}
+O\left(\frac{1}{s^{2(l+1)+\mu}}\right)\right\}\\
\\
\displaystyle
=i(1-\zeta)^2
\left\{
\sum_{n=0}^l\frac{(k^2e^{2i\theta})^n(n+1)(n+1+\mu)}{s^{2n}}
+O\left(\frac{1}{s^{2(l+1)}}\right)\right\}
\frac{ie^{i(\theta+\frac{\pi}{2})\mu}k^{\mu}}{s^{\mu}}.
\end{array}
\tag {2.23}
$$

Let $\vert\zeta\vert<1$.
Since
$$\displaystyle
\sum_{n=0}^{\infty}\zeta^n(n+1)=(1-\zeta)^{-2},\,\,
\sum_{n=0}^{\infty}\zeta^n(n+1)^2
=(1+\zeta)(1-\zeta)^{-3},
$$
we have
$$\begin{array}{c}
\displaystyle
\sum_{n=0}^\infty
\zeta^n(n+1)(n+1+\mu)=\{1+\zeta+\mu(1-\zeta)\}(1-\zeta)^{-3}
\end{array}
$$
and for each fixed $l$
$$\begin{array}{c}
\displaystyle
\sum_{n=0}^l
\zeta^n(n+1)(n+1+\mu)=\{1+\zeta+\mu(1-\zeta)\}(1-\zeta)^{-3}
+O(\vert\zeta\vert^{l+1}).
\end{array}
\tag {2.24}
$$
Substituting (2.24) with $\zeta=(k/s)^2e^{2i\theta}$ into (2.23), we have
$$\begin{array}{c}
\displaystyle
(\tau\,\cos\,\theta-i\sqrt{\tau^2+k^2}\,\sin\,\theta)^2K_{\mu}(\tau,\theta)\\
\\
\displaystyle
=i(1-\zeta)^2
\{1+\zeta+\mu(1-\zeta)\}(1-\zeta)^{-3}
\frac{ie^{i(\theta+\frac{\pi}{2})\mu}k^{\mu}}{s^{\mu}}
+O\left(\frac{1}{s^{2(l+1)+\mu}}\right).
\end{array}
$$
Since $l$ can be arbitrary large whenever $n$ is fixed,
$\mu=\mu_m$ and $0\le m\le n$, we obtain (2.16).

\noindent
$\Box$

We continue the proof of (2.3).
One can write
$$
\displaystyle
\tau\,\sin\,\theta+i\sqrt{\tau^2+k^2}\,\cos\,\theta
=\frac{ise^{-i\theta}}{2}\left\{1+\left(\frac{k}{s}\right)^2e^{2i\theta}\right\}
$$
and
$$
\displaystyle
\tau\cos\,\theta-i\sqrt{\tau^2+k^2}\,\sin\,\theta
=\frac{se^{-i\theta}}{2}
\left\{1-\left(\frac{k}{s}\right)^2 e^{2i\theta}\right\}.
$$
From these and (2.15) one gets
$$\begin{array}{c}
\displaystyle
\left\{(\tau\,\sin\,\theta+i\sqrt{\tau^2+k^2}\,\cos\,\theta)
-x_0\cdot c_{\tau}(\omega^{\perp})
(\tau\,\cos\,\theta-i\sqrt{\tau^2+k^2}\,\sin\,\theta)\right\}
I_{\mu}(\tau,\theta)\\
\\
\displaystyle
=\left\{i(1+\zeta)
(1-\zeta)^{-1}-x_0\cdot c_{\tau}(\omega^{\perp})\right\}
\frac{ie^{i(\theta+\frac{\pi}{2})\mu}k^{\mu}}{s^{\mu}}
+O(s^{-\infty}).
\end{array}
\tag {2.25}
$$
Since
$$
\displaystyle
i(1+\zeta)(1-\zeta)^{-1}
-i(1-\zeta)^2\{1+\zeta+\mu(1-\zeta)\}(1-\zeta)^{-3}
=-i\mu,
$$
it follows from (2.13), (2.14), (2.16) and (2.25) that
$$\begin{array}{c}
\displaystyle
\sqrt{\tau^2+k^2}I_p(\tau)
=\sum_{m=2}^n
\alpha_m(i\mu_m+x_0\cdot c_{\tau}(\omega^{\perp}))
\frac{e^{i(p+\frac{\pi}{2})\mu_m}k^{\mu_m}}{s^{\mu_m}}
+O\left(\frac{1}{s^{\mu_{n+1}}}\right)
\end{array}
$$
and
$$\begin{array}{c}
\displaystyle
\sqrt{\tau^2+k^2}I_q(\tau)
=\sum_{m=2}^n
\alpha_m(-1)^m(i\mu_m+x_0\cdot c_{\tau}(\omega^{\perp}))
\frac{e^{i(q+\frac{\pi}{2})\mu_m}k^{\mu_m}}{s^{\mu_m}}
+O\left(\frac{1}{s^{\mu_{n+1}}}\right).
\end{array}
$$
Now from these and (2.6) we obtain (2.3).

\noindent
$\Box$

\section{Proof of Theorem 1.2}

First consider the case when every end points of $\Sigma_1$, $\Sigma_2$, $\cdots$, $\Sigma_m$
satisfies $x\cdot\omega<h_{\Sigma}(\omega)$.  Then $x_0\in\Sigma$ with $x_0\cdot\omega=h_{\Sigma}(\omega)$
should be a vertex of $D$ and a point where two segments in some $\Sigma_j$ meet.
We take the same polar coordinates as those of Section 2.

Integration by parts gives
$$\displaystyle
I_{\Sigma}'(\tau;\omega,d,k)
=\int_{\Sigma}[u]\frac{\partial}{\partial\nu}\partial_{\tau}vdS,
$$
where $[u]=u^+\vert_{\partial D}-u^-\vert_{\partial D}$.
Localizing this integral at $x_0$, we have, modulo exponentially decaying as
$\tau\longrightarrow\infty$
$$\begin{array}{c}
\displaystyle
e^{-i\sqrt{\tau^2+k^2}\,x_0\cdot\omega^{\perp}}e^{-\tau h_D(\omega)}I_{\Sigma}'(\tau;\omega,d,k)\\
\\
\displaystyle
\sim
e^{-i\sqrt{\tau^2+k^2}\,x_0\cdot\omega^{\perp}}e^{-\tau h_D(\omega)}
\int_{\Gamma_p}
[u]\frac{\partial}{\partial\nu}\partial_{\tau}vdS
+e^{-i\sqrt{\tau^2+k^2}\,x_0\cdot\omega^{\perp}}e^{-\tau h_D(\omega)}\int_{\Gamma_q}
[u]\frac{\partial}{\partial\nu}\partial_{\tau}vdS
\\
\\
\displaystyle
\equiv
I_p(\tau)+I_q(\tau).
\end{array}
\tag {3.1}
$$

It follows from (1.3) and (2.7) that
$$\begin{array}{c}
\displaystyle
\sqrt{\tau^2+k^2}I_p(\tau)
=
i\left\{(\tau\sin\,p+i\sqrt{\tau^2+k^2}\cos\,p)
-x_0\cdot c_{\tau}(\omega^{\perp})(\tau\cos\,p-i\sqrt{\tau^2+k^2}\sin\,p)\right\}\\
\\
\displaystyle
\times
\int_0^{\eta}(u^+(r,0)-u^-(r,2\pi))e^{\tau r\,\sin p}e^{i\sqrt{\tau^2+k^2}\,r\,\cos\,p}dr\\
\\
\displaystyle
-i(\tau\cos\,p-i\sqrt{\tau^2+k^2}\sin\,p)^2
\int_0^{\eta}r(u^+(r,0)-u^-(r,2\pi))e^{\tau r\,\sin p}e^{i\sqrt{\tau^2+k^2}\,r\,\cos\,p}dr.
\end{array}
\tag {3.2}
$$

It follows also from (1.3) and (2.8) that
$$\begin{array}{c}
\displaystyle
\sqrt{\tau^2+k^2}I_q(\tau)
=
-i\left\{(\tau\sin\,q+i\sqrt{\tau^2+k^2}\cos\,q)
-x_0\cdot c_{\tau}(\omega^{\perp})(\tau\cos\,q-i\sqrt{\tau^2+k^2}\sin\,q)\right\}\\
\\
\displaystyle
\times
\int_0^{\eta}(u^+(r,\Theta)-u^-(r,\Theta))e^{\tau r\,\sin q}e^{i\sqrt{\tau^2+k^2}\,r\,\cos\,q}dr\\
\\
\displaystyle
+i(\tau\cos\,q-i\sqrt{\tau^2+k^2}\sin\,q)^2
\int_0^{\eta}r(u^+(r,\Theta)-u^-(r,\Theta))e^{\tau r\,\sin q}e^{i\sqrt{\tau^2+k^2}\,r\,\cos\,q}dr.
\end{array}
\tag {3.3}
$$

By Proposition 4.4 in \cite{IE4}, we have
$$\begin{array}{c}
\displaystyle
u^+(r,\theta)
=\sum_{m=1}^{\infty}\alpha_m^+ J_{\mu_m^+}(kr)\cos\,\mu_m^+\theta, 0<r<\eta, 0\le \theta\le \Theta,\\
\\
\displaystyle
u^-(r,\theta)
=\sum_{m=1}^{\infty}\alpha_m^- J_{\mu_m^-}(kr)\cos\,\mu_m^-(\theta-\Theta), 0<r<\eta, \Theta\le \theta\le 2\pi,
\end{array}
$$
where $\mu_m^+=(m-1)\pi/\Theta$ and
$\mu_m^-=(m-1)\pi/(2\pi-\Theta)$. For the precise meaning of this
expansion see \cite{IE4}. From these and $\mu_m^+<\mu_m^-$ we have
$$\begin{array}{c}
\displaystyle
[u]_p
=u^+(r,0)-u^-(r,2\pi)
=\sum_{m=1}^l\{J_{\mu_m^+}(kr)\alpha_m^+
+J_{\mu_m^-}(kr)\alpha_m^-(-1)^m\}+O(r^{\mu_{l+1}^+}),
\\
\\
\displaystyle
[u]_q
=u^+(r,\Theta)-u^-(r,\Theta)
=-\sum_{m=1}^l\{J_{\mu_m^+}(kr)\alpha_m^+(-1)^m
+J_{\mu_m^-}(kr)\alpha_m^-\}+O(r^{\mu_{l+1}^+}).
\end{array}
\tag {3.4}
$$
It follows from (3.4) and (2.11) that
$$\begin{array}{c}
\displaystyle
\int_0^{\eta}(u^+(r,0)-u^-(r,2\pi))e^{\tau r\,\sin p}e^{i\sqrt{\tau^2+k^2}\,r\,\cos\,p}dr\\
\\
\displaystyle
=\sum_{m=1}^n\{\alpha_m^+I_{\mu_m^+}(\tau,p)
+\alpha_m^-(-1)^mI_{\mu_m^-}(\tau,p)\}
+O(\tau^{-(\mu_{n+1}^++1)}),
\end{array}
$$

$$\begin{array}{c}
\displaystyle
\int_0^{\eta}r(u^+(r,0)-u^-(r,2\pi))e^{\tau r\,\sin p}e^{i\sqrt{\tau^2+k^2}\,r\,\cos\,p}dr\\
\\
\displaystyle
=\sum_{m=1}^n\{\alpha_m^+K_{\mu_m^+}(\tau,p)
+\alpha_m^-(-1)^mK_{\mu_m^-}(\tau,p)\}
+O(\tau^{-(\mu_{n+1}^++2)}),
\end{array}
$$

$$\begin{array}{c}
\displaystyle
\int_0^{\eta}(u^+(r,\Theta)-u^-(r,\Theta))e^{\tau r\,\sin q}e^{i\sqrt{\tau^2+k^2}\,r\,\cos\,q}dr\\
\\
\displaystyle
=-\sum_{m=1}^n\{\alpha_m^+(-1)^mI_{\mu_m^+}(\tau,q)
+\alpha_m^-I_{\mu_m^-}(\tau,q)\}
+O(\tau^{-(\mu_{n+1}^++1)}),
\end{array}
$$

$$\begin{array}{c}
\displaystyle
\int_0^{\eta}r(u^+(r,\Theta)-u^-(r,\Theta))e^{\tau r\,\sin q}e^{i\sqrt{\tau^2+k^2}\,r\,\cos\,q}dr\\
\\
\displaystyle
=-\sum_{m=1}^n\{\alpha_m^+(-1)^mK_{\mu_m^+}(\tau,q)
+\alpha_m^-K_{\mu_m^-}(\tau,q)\}
+O(\tau^{-(\mu_{n+1}^++2)}).
\end{array}
$$

Substituting these into (3.2) and (3.3), we obtain

$$\begin{array}{c}
\displaystyle
\sqrt{\tau^2+k^2}I_p(\tau)\\
\\
\displaystyle
=
i\left\{(\tau\sin\,p+i\sqrt{\tau^2+k^2}\cos\,p)
-x_0\cdot c_{\tau}(\omega^{\perp})(\tau\cos\,p-i\sqrt{\tau^2+k^2}\sin\,p)\right\}\\
\\
\displaystyle
\times
\sum_{m=1}^n\{\alpha_m^+I_{\mu_m^+}(\tau,p)
+\alpha_m^-(-1)^mI_{\mu_m^-}(\tau,p)\}
\\
\\
\displaystyle
-i(\tau\cos\,p-i\sqrt{\tau^2+k^2}\sin\,p)^2
\sum_{m=1}^n\{\alpha_m^+K_{\mu_m^+}(\tau,p)
+\alpha_m^-(-1)^mK_{\mu_m^-}(\tau,p)\}
+O(\tau^{-\mu_{n+1}^+})
\end{array}
\tag {3.5}
$$
and
$$\begin{array}{c}
\displaystyle
\sqrt{\tau^2+k^2}I_q(\tau)\\
\\
\displaystyle
=
-i\left\{(\tau\sin\,q+i\sqrt{\tau^2+k^2}\cos\,q)
-x_0\cdot c_{\tau}(\omega^{\perp})(\tau\cos\,q-i\sqrt{\tau^2+k^2}\sin\,q)\right\}\\
\\
\displaystyle
\times
\sum_{m=1}^n\{\alpha_m^+(-1)^mI_{\mu_m^+}(\tau,q)
+\alpha_m^-I_{\mu_m^-}(\tau,q)\}
\\
\\
\displaystyle
+i(\tau\cos\,q-i\sqrt{\tau^2+k^2}\sin\,q)^2
\sum_{m=1}^n\{\alpha_m^+(-1)^mK_{\mu_m^+}(\tau,q)
+\alpha_m^-K_{\mu_m^-}(\tau,q)\}
+O(\tau^{-\mu_{n+1}^+}).
\end{array}
\tag {3.6}
$$

It follows from (3.5), (3.6), (2.16) and (2.25) that
$$\begin{array}{c}
\displaystyle
\sqrt{\tau^2+k^2}I_p(\tau)
=\sum_{m=1}^n
\alpha_m^+(i\mu_m^++x_0\cdot c_{\tau}(\omega^{\perp}))
\frac{e^{i(p+\frac{\pi}{2})\mu_m^+}k^{\mu_m^+}}{s^{\mu_m^+}}\\
\\
\displaystyle
+\sum_{m=1}^n
\alpha_m^-(-1)^m(i\mu_m^-+x_0\cdot c_{\tau}(\omega^{\perp}))
\frac{e^{i(p+\frac{\pi}{2})\mu_m^-}k^{\mu_m^-}}{s^{\mu_m^-}}
+O\left(\frac{1}{s^{\mu_{n+1}^+}}\right)
\end{array}
$$
and
$$\begin{array}{c}
\displaystyle
\sqrt{\tau^2+k^2}I_q(\tau)
=\sum_{m=1}^n
\alpha_m^+(-1)^m(i\mu_m^++x_0\cdot c_{\tau}(\omega^{\perp}))
\frac{e^{i(q+\frac{\pi}{2})\mu_m^+}k^{\mu_m^+}}{s^{\mu_m^+}}\\
\\
\displaystyle
+\sum_{m=1}^n
\alpha_m^-(i\mu_m^-+x_0\cdot c_{\tau}(\omega^{\perp}))
\frac{e^{i(q+\frac{\pi}{2})\mu_m^-}k^{\mu_m^-}}{s^{\mu_m^-}}
+O\left(\frac{1}{s^{\mu_{n+1}^+}}\right).
\end{array}
$$
Now it follows from this and (3.1) that
$$\begin{array}{c}
\displaystyle
\sqrt{\tau^2+k^2}e^{-i\sqrt{\tau^2+k^2}\,x_0\cdot\omega^{\perp}}e^{-\tau h_D(\omega)}I_{\Sigma}'(\tau;\omega,d,k)\\
\\
\displaystyle
=-i\sum_{m=1}^n
\alpha_m^+\{e^{ip\mu_m^+}+(-1)^me^{iq\mu_m^+}\}(-\mu_m^++ix_0\cdot c_{\tau}(\omega^{\perp}))
\frac{e^{i\frac{\pi}{2}\mu_m^+}k^{\mu_m^+}}{s^{\mu_m^+}}\\
\\
\displaystyle
-i\sum_{m=1}^n
\alpha_m^-\{(-1)^me^{ip\mu_m^-}+e^{iq\mu_m^-}\}(-\mu_m^-+ix_0\cdot c_{\tau}(\omega^{\perp}))
\frac{e^{i\frac{\pi}{2}\mu_m^-}k^{\mu_m^-}}{s^{\mu_m^-}}
+O\left(\frac{1}{s^{\mu_{n+1}^+}}\right).
\end{array}
\tag {3.7}
$$
Since we have the following cancellation (\cite{IE4})
$$\displaystyle
(-1)^me^{ip\mu_m^-}+e^{iq\mu_m^-}=0,
\tag {3.8}
$$
(3.7) becomes
$$\begin{array}{c}
\displaystyle
\sqrt{\tau^2+k^2}e^{-i\sqrt{\tau^2+k^2}\,x_0\cdot\omega^{\perp}}e^{-\tau h_D(\omega)}I_{\Sigma}'(\tau;\omega,d,k)\\
\\
\displaystyle
=-i\sum_{m=1}^n
\alpha_m^+\{e^{ip\mu_m^+}+(-1)^me^{iq\mu_m^+}\}(-\mu_m^++ix_0\cdot c_{\tau}(\omega^{\perp}))
\frac{e^{i\frac{\pi}{2}\mu_m^+}k^{\mu_m^+}}{s^{\mu_m^+}}
+O\left(\frac{1}{s^{\mu_{n+1}^+}}\right).
\end{array}
\tag {3.9}
$$

By the way, we have already known in \cite{IE4} that
$$\begin{array}{c}
\displaystyle
e^{-i\sqrt{\tau^2+k^2}\,x_0\cdot\omega^{\perp}}e^{-\tau h_D(\omega)}I_{\Sigma}(\tau;\omega,d,k)\\
\\
\displaystyle
=-i\sum_{m=1}^n
\alpha_m^+\{e^{ip\mu_m^+}+(-1)^me^{iq\mu_m^+}\}
\frac{e^{i\frac{\pi}{2}\mu_m^+}k^{\mu_m^+}}{s^{\mu_m^+}}
+O\left(\frac{1}{s^{\mu_{n+1}^+}}\right)
\end{array}
\tag {3.10}
$$
and there exists a $m\ge 2$ such that
$\alpha_m^+\{e^{ip\mu_m^+}+(-1)^me^{iq\mu_m^+}\}\not=0$.
Having these together with (3.9), hereafter
we take the same course as the obstacle case.

The case when there is an end point of some $\Sigma_j$ such that
$x\cdot\omega=h_{\Sigma}(\omega)$ corresponds to the case when
$p=q$.  We omit its description.

\noindent
$\Box$

{\bf\noindent Remark 3.1.}  From (3.9) and (3.10) we see that the field $u^-(r,\theta)$ never affects
the asymptotic behaviour of $I_{\Sigma}'(\tau;\omega,d,k)$ and
$I_{\Sigma}(\tau;\omega,d,k)$ as $\tau\longrightarrow\infty$ modulo {\it rapidly decreasing}.
The key point is cancellation (3.8).

\section{Some other applications}

In this section we present some implications of the argument done for the proof of Theorems 1.1.

\subsection{From the far-field pattern of the scattered wave for a single incident plane wave.}
Let $u=e^{ikx\cdot d}+w$ be the same as that of Theorem 1.1.
It is well known that $w$ has the asymptotic expansion
as $r\longrightarrow\infty$ uniformly with respect to $\varphi\in\,S^1$:
$$\displaystyle
w(r\varphi)=\frac{e^{ikr}}{\sqrt{r}}F(\varphi;d,k)+O\left(\frac{1}{r^{3/2}}\right).
$$
The coefficient $F(\varphi;d,k)$ is called the far-field pattern of
the scattered wave $w$ at direction $\varphi$.

In this subsection we present a direct formula that extracts the
coordinates of the vertices of the convex hull of unknown
polygonal sound hard obstacles $D=D_1\cup\cdots D_m$ from the far field pattern for
fixed $d$ and $k$.

We identify $\varphi=(\varphi_1,\varphi_2)$ with the complex number
given by $\varphi_1+i\varphi_2$ and denote it by the same symbol $\varphi$.
Given $N=1,\cdots$, $\tau>0$, $\omega\in S^1$ and $k>0$ define the function $g_N(\,\cdot\,;\tau,k,\omega)$ on $S^1$
by the formula
$$
\displaystyle
g_{N}(\varphi;\tau,k,\omega)
=\frac{1}{2\pi}
\sum_{\vert m\vert\le N}
\left\{\frac{ik\varphi}{(\tau+\sqrt{\tau^2+k^2})\omega}\right\}^m.
$$
Then we have
$$\displaystyle
\partial_{\tau}g_{N}(\varphi;\tau,k,\omega)
=-
\frac{1}{2\pi\sqrt{\tau^2+k^2}}
\sum_{1\le\vert m\vert\le N}
m\left\{\frac{ik\varphi}{(\tau+\sqrt{\tau^2+k^2})\omega}\right\}^m.
$$
In this subsection $B_R$ denotes the open disc centered at the
origin of the coordinates with radius $R$ and we assume that
$\overline D\subset B_R$.

\proclaim{\noindent Theorem 4.1.}  Let $\omega$ be regular with
respect to $D$.  Let $\beta_0$ be the unique positive solution of the equation
$$\displaystyle
\frac{2}{e}s+\log s=0.
$$
Let $\beta$ satisfy $0<\beta<\beta_0$.  Let $\{\tau(N)\}_{N=1,\cdots}$ be an arbitrary
sequence of positive numbers satisfying, as $N\longrightarrow\infty$
$$\displaystyle
\tau(N)=\frac{\beta N}{eR}+O(1).
$$
Then the formula
$$\displaystyle
\lim_{N\longrightarrow\infty}
\frac
{\displaystyle\int_{S^1}F(-\varphi;d,k)\partial_{\tau}g_N(\varphi;\tau(N),k,\omega)dS(\varphi)}
{\displaystyle\int_{S^1}F(-\varphi;d,k)g_N(\varphi;\tau(N),k,\omega)dS(\varphi)}
=h_D(\omega)+ix_0\cdot\omega^{\perp},
\tag {4.1}
$$
is valid.
\endproclaim

In \cite{Ik11} we have shown that, under the same choice of $\tau(N)$ and
$\omega$ being regular with respec to $D$
$$\displaystyle
\lim_{N\longrightarrow\infty}
\frac{1}{\tau(N)}
\log\left\vert\int_{S^1}F(-\varphi;d,k)g_N(\varphi;\tau(N),k,\omega)dS(\varphi)\right\vert
=h_D(\omega).
$$
Thus Theorem 4.1 corresponds to Theorem 1.1.

Let us describe the proof of Theorem 4.1.
The starting point is the following identity which is a consequence of the formula (2.9)
in \cite{CK1}:
$$\begin{array}{c}
\displaystyle
-\frac{\sqrt{8\pi k}}{e^{i\pi/4}}
\int_{S^1}F(-\varphi;d,k)g_N(\varphi;\tau,k,\omega)dS(\varphi)\\
\\
\displaystyle
=I(\tau;\omega,d,k)
+\int_{\partial B_R}\left\{\frac{\partial u}{\partial\nu}(v_{g_N}-v_{\tau})
-\frac{\partial}{\partial\nu}(v_{g_N}-v_{\tau})u\right\}dS
\end{array}
\tag {4.2}
$$
where
$$\displaystyle
v_{g_N}(y)=\int_{S^1}e^{iy\cdot\varphi}g_N(\varphi;\tau,k,\omega)dS(\varphi).
$$
By Theorem 2.1 in \cite{Ik11} we know that the second term in the right-hand side of (4.2) has the bound $e^{-\tau(N)R}O(N^{-\infty})$
as $N\longrightarrow\infty$.  Since $e^{-\tau(N)h_D(\omega)}=O(e^{\tau(N)R})$, it follows from (2.4) and (4.2) that
$$\begin{array}{c}
\displaystyle
\lim_{N\longrightarrow\infty}\tau(N)^{\mu_{m^*}}e^{-i\sqrt{\tau(N)^2+k^2}x_0\cdot\omega^{\perp}}e^{-\tau(N) h_D(\omega)}\frac{\sqrt{8\pi k}}{e^{i\pi/4}}
\int_{S^1}F(-\varphi;d,k)g_N(\varphi;\tau(N),k,\omega)dS(\varphi)\\
\\
\displaystyle
=i\beta e^{i\frac{\pi}{2}\mu_{m^*}}\left(\frac{k}{2}\right)^{\mu_{m^*}}.
\end{array}
\tag {4.3}
$$
Similarly to (4.2) we have
$$\begin{array}{c}
\displaystyle
-\frac{\sqrt{8\pi k}}{e^{i\pi/4}}
\int_{S^1}F(-\varphi;d,k)\partial_{\tau}g_N(\varphi;\tau,k,\omega)dS(\varphi)\\
\\
\displaystyle
=I'(\tau;\omega,d,k)
+\int_{\partial B_R}\left\{\frac{\partial u}{\partial\nu}(v_{\partial_{\tau}g_N}-\partial_{\tau}v_{\tau})
-\frac{\partial}{\partial\nu}(v_{\partial_{\tau}g_N}-\partial_{\tau}v_{\tau})u\right\}dS
\end{array}
\tag {4.4}
$$
where
$$\displaystyle
v_{\partial_{\tau}g_N}(y)
=\int_{S^1}e^{iy\cdot\varphi}\partial_{\tau}g_N(\varphi;\tau, \omega,k,d)dS(\varphi).
$$

The following lemma corresponds to Theorem 2.1 in \cite{Ik11} and
see Appendix for the proof.

\proclaim{\noindent Lemma 4.1.}
We have, as $N\longrightarrow\infty$
$$\begin{array}{c}
\displaystyle
e^{R\tau(N)}\sup_{\vert y\vert\le R}\left\vert\int_{S^1}e^{iky\cdot\varphi}\partial_{\tau}g_N(\varphi;\tau(N),k,\omega)dS(\varphi)
-\partial_{\tau}v_{\tau}(y)\vert_{\tau=\tau(N)}\right\vert
\\
\\
\displaystyle
+e^{R\tau(N)}\sup_{\vert y\vert\le R}\left\vert\nabla
\left\{\int_{S^1}e^{iky\cdot\varphi}\partial_{\tau}g_N(\varphi;\tau(N),k,\omega)dS(\varphi)
-\partial_{\tau}v_{\tau}(y)\vert_{\tau=\tau(N)}\right\}\right\vert
\\
\\
\displaystyle
=O(N^{-\infty}).
\end{array}
$$
\endproclaim

Now from Lemma 4.1, (4.3) and (2.5) we obtain
$$\begin{array}{c}
\displaystyle
\lim_{N\longrightarrow\infty}\tau(N)^{\mu_{m^*}}e^{-i\sqrt{\tau(N)^2+k^2}x_0
\cdot\omega^{\perp}}e^{-\tau(N) h_D(\omega)}\frac{\sqrt{8\pi k}}{e^{i\pi/4}}
\int_{S^1}F(-\varphi;d,k)\partial_{\tau}g_N(\varphi;\tau(N),k,\omega)dS(\varphi)\\
\\
\displaystyle
=i\beta(x_0\cdot\omega+ix_0\cdot\omega^{\perp})e^{i\frac{\pi}{2}\mu_{m^*}}
\left(\frac{k}{2}\right)^{\mu_{m^*}}.
\end{array}
$$
From this together with (4.4) yields (4.1).

\subsection{From the Cauchy data of the scattered wave for a single point source}

Let $y\in\Bbb R^2\setminus\overline D$.   Let $E=E_D(x,y)$ be the unique solution of the scattering problem:
$$\begin{array}{c}
\displaystyle
(\triangle+k^2)E=0\,\,\text{in}\,\Bbb R^2\setminus\overline D,\\
\\
\displaystyle
\frac{\partial}{\partial\nu}E=-\frac{\partial}{\partial\nu}\Phi_0(\,\cdot\,,y)
\,\text{on}\,\partial D,\\
\\
\displaystyle
\lim_{r\longrightarrow\infty}\sqrt{r}\left(\frac{\partial E}{\partial r}-ikE\right)=0,
\end{array}
$$
where
$$\displaystyle
\Phi_0(x,y)=\frac{i}{4}H^{(1)}_0(k\vert x-y\vert)
$$
and $H^{(1)}_0$ denotes the Hankel function of the first kind \cite{O}.

The total wave outside $D$ exerted by the point source located at $y$ is given by the formula:
$$\displaystyle
\Phi_D(x,y)=\Phi_0(x,y)+E_D(x,y),\,x\in\Bbb R^2\setminus\overline D.
$$

In this subsection we consider the following problem.

\noindent
{\bf Inverse Problem.}  Let $R_1>R$.
We denote by $B_R$ and $B_{R_1}$ the open discs centered at a common point with radius $R$ and $R_1$, respectively.
Assume that $\overline D\subset B_R$.
Fix $k>0$ and $y\in\partial B_{R_1}$.
Extract information about the location and shape
of $D$ from $\Phi_D(x,y)$ given at all $x\in\partial B_{R}$.

Define
$$\displaystyle
J(\tau; \omega, y, k)
=
\int_{\partial B_{R}}
\left(
\frac{\partial}{\partial\nu}\Phi_D(x,y)\cdot v_{\tau}(x;\omega)
-\frac{\partial}{\partial\nu}v_{\tau}(x;\omega)\cdot\Phi_D(x,y)\right)dS(x).
$$
Then we have
$$\displaystyle
J'(\tau; \omega, y, k)
=
\int_{\partial B_{R}}
\left(
\frac{\partial}{\partial\nu}\Phi_D(x,y)\cdot \partial_{\tau}v_{\tau}
-\frac{\partial}{\partial\nu}(\partial_{\tau}v_{\tau})\cdot\Phi_D(x,y)\right)dS(x).
$$
Note that $(\partial/\partial\nu)\Phi_D(x,y)$ for $x\in\partial B_R$ can be computed
from $\Phi_D(x,y)$ given at all $x\in\partial B_R$ by solving an exterior Dirichlet problem 
for the Helmholtz equation.  See \cite{IP} for this point.
This remark applies also to $(\partial/\partial\nu)u$ on $\partial B_R$ in Theorems 1.1
and 1.2.

A combination of the proof of Theorem 1.2 in \cite{IP} and the same argument as done in the proof of Theorem 1.1 yields the following formula.

\proclaim{\noindent Theorem 4.2.}
Assume that
$$\displaystyle
\text{diam}\,D<\text{dist}\,(D,\partial B_{R_1}).
\tag {4.5}
$$
Let $\omega$ be regular with respect to $D$.
Let $x_0\in\partial D$ be the point with $x_0\cdot\omega=h_D(\omega)$.
Then, there exists a $\tau_0>0$ such that, for all $\tau\ge\tau_0$
$\vert J(\tau;\omega,d,k)\vert>0$ and the formula
$$\displaystyle
\lim_{\tau\longrightarrow\infty}
\frac{J'(\tau;\omega,y,k)}
{J(\tau;\omega,y,k)}
=h_D(\omega)+ix_0\cdot\omega^{\perp},
$$
is valid.
\endproclaim

Note that it is an open problem whether one can drop condition (4.5).  For more information about this see \cite{IP}.

\centerline{{\bf Acknowledgements}}

This research was partially supported by Grant-in-Aid for
Scientific Research (C)(No. 21540162) of Japan  Society for the
Promotion of Science.
The author thanks Takashi Ohe for useful discussion.

\section{Appendix}

\subsection{Proof of Proposition 2.1}

{\it\noindent Proof.}
For each $m$
define
$$\displaystyle
u_m(r,\theta)=\int_0^{2\Theta}u(r,\theta)\cos\,\mu_m\theta d\theta,\,\,(r,\theta)
\in\,]0,\,2\eta[\,\times
\,[0,\,\Theta].
$$
Then we see that $u_m$ satisfies the equation
$$\displaystyle
-(ry')'+(\frac{\mu_m^2}{r}-k^2r)y=0\,\,\text{in}\,]0,\,2\eta[
$$
and thus this yields that there exist numbers $\alpha_m$, $\beta_m$ such that
$$
\displaystyle
u_m(r)=\alpha_m J_{\mu_m}(kr)+\beta_m Y_{\mu_m}(kr),
$$
where $Y_{\mu_m}(kr)$ denotes the Bessel function of the second kind (\cite{O}).
Then a similar argument done in \cite{G} and the behaviour
of $J_{\mu_m}(kr)$ and $Y_{\mu_m}(kr)$ as
$r\longrightarrow 0$ (\cite{O}), one concludes that $\beta_m=0$ for all $m\ge
1$.  Substituting $u_m(r)=\alpha_m J_{\mu_m}(kr)$ into the
inequality
$$\displaystyle
\int_0^{2\eta}r\vert u'_m(r)\vert^2 dr\le\Vert\nabla u\Vert^2_{L^2(B_{2\eta}(x_0)\cap
(B_R\setminus\overline D))}<\infty,
$$
we obtain
$$\displaystyle
\alpha_m^2\int_0^{2\eta}r\vert\{J_{\mu_m}(kr)\}'\vert^2dr\le C
\tag {A.1}
$$
where $C$ is a positive constant independent of $m$.

Thus the problem is to estimate the integral in the left-hand side
of (A.1) from below. For the purpose we make use of the following
formula which can be checked directly:
$$\displaystyle
\{J_{\mu_m}(kr)\}'
=\frac{\mu_m}{r}J_{\mu_m}(kr)-kJ_{\mu_{m}+1}(kr).
$$
From this we have
$$
\displaystyle
r\vert\{J_{\mu_m}(kr)\}'\vert^2
=\frac{\mu_m^2}{r}
\vert J_{\mu_m}(kr)\vert^2
-2\mu_mkJ_{\mu_m}(kr)J_{\mu_m+1}(kr)
+k^2r\vert J_{\mu_m+1}(kr)\vert^2
$$
and thus
$$\displaystyle
\int_0^{2\eta}r\vert\{J_{\mu_m}(kr)\}'\vert^2dr
\ge\mu_m^2\int_0^{2\eta}\frac{1}{r}
\vert J_{\mu_m}(kr)\vert^2dr
-2\mu_mk\int_o^{2\eta}J_{\mu_m}(k)J_{\mu_m+1}(kr)dr.
\tag {A.2}
$$
By formula (37) on p.338 in \cite{B} and a change of independent variable
we have
$$\begin{array}{c}
\displaystyle
k\int_0^{2\eta}J_{\mu_m}(kr)J_{\mu_m+1}(kr)dr
=\int_0^{2\eta k}J_{\mu_m}(r)J_{\mu_m+1}(r)dr\\
\\
\displaystyle
=\sum_{n=0}^{\infty}\vert J_{\mu_m+n+1}(2\eta k)\vert^2.
\end{array}
\tag {A.3}
$$
Furthermore from (A.3) and the recurrence relation of the Bessel
functions
$$\displaystyle
\mu_mJ_{\mu_m}(kr)=kr(J_{\mu_m-1}(kr)+J_{\mu_{m}+1}(kr))
$$
which can be checked directly
we have
$$\begin{array}{c}
\displaystyle
\int_0^{2\eta}\frac{1}{r}\vert J_{\mu_m}(kr)\vert^2 dr
=\frac{k}{\mu_m}\int_0^{2\eta}J_{\mu_m-1}(kr)J_{\mu_m}(kr)dr
+\frac{k}{\mu_m}\int_0^{2\eta}J_{\mu_m}(kr)J_{\mu_m+1}(kr)dr\\
\\
\displaystyle
=\frac{1}{\mu_m}\int_0^{2\eta k}J_{\mu_m-1}(r)J_{\mu_m}(r)dr
+\frac{1}{\mu_m}\int_0^{2\eta k}J_{\mu_m}(r)J_{\mu_m+1}(r)dr\\
\\
\displaystyle
=\frac{1}{\mu_m}\sum_{n=0}^{\infty}\vert J_{\mu_m+n}(2\eta k)\vert^2
+\frac{1}{\mu_m}\sum_{n=0}^{\infty}\vert J_{\mu_m+n+1}(2\eta k)\vert^2\\
\\
\displaystyle
=
\frac{1}{\mu_m}\vert J_{\mu_m}(2\eta k)\vert^2
+\frac{2}{\mu_m}\sum_{n=0}^{\infty}\vert J_{\mu_m+n+1}(2\eta k)\vert^2.
\end{array}
$$
This together with (A.2) and (A.3) yields
$$
\displaystyle
\int_0^{2\eta}r\vert \{J_{\mu_m}(kr)\}'\vert^2dr
\ge \mu_m\vert J_{\mu_m}(2\eta k)\vert^2.
\tag {A.4}
$$
Since
$$\displaystyle
J_{\mu_m}(2\eta k)
=\frac{(\eta k)^{\mu_m}}{\Gamma(1+\mu_m)}
\left\{1+\sum_{n=1}^{\infty}
\frac{(-1)^m\Gamma(1+\mu_m)}{\Gamma(1+\mu_m+n)}(\eta k)^{2n}\right\}
$$
and
$$\displaystyle
\frac{\Gamma(1+\mu_m)}
{\Gamma(1+\mu_m+n)}
=\frac{1}{\Pi_{j=1}^n(j+\mu_m)}
\le\frac{1}{n!},
$$
it holds that
$$\displaystyle
\vert J_{\mu_m}(2\eta k)\vert
\ge\frac{(\eta k)^{\mu_m}}
{\Gamma(1+\mu_m)}
(1-(e^{\eta k}-1)).
$$
Thus choosing $\eta k$ in such a way that $e^{\eta k}-1<1/2$, that is, $\eta k<\log(3/2)$, we obtain
$$\displaystyle
\vert J_{\mu_m}(2\eta k)\vert
\ge\frac{1}{2}\frac{(\eta k)^{\mu_m}}
{\Gamma(1+\mu_m)}.
$$
This together with (A.4) gives the following estimate:
$$\displaystyle
\int_0^{2\eta}r\vert\{J_{\mu_m}(kr)\}'\vert^2 dr
\ge\frac{1}{4}\frac{\mu_m}{\Gamma(1+\mu_m)^2}(\eta k)^{2\mu_m}.
\tag {A.5}
$$
Then (A.1) and (A.5) give the estimate in (2).
The remaining parts of the statements are
consequences of (2) and the completeness of $(\cos\,\mu_m\theta)_{m=1}^{\infty}$
in $L^2(]0,\,\Theta[)$.

\noindent
$\Box$

\subsection{Proof of Lemma 2.3}
Since
$$
\displaystyle
\Gamma(n+2+l+\mu)
=\{\Pi_{j=1}^{n+1}(j+l+\mu)\}\Gamma(1+l+\mu),
$$
one can rewrite
$$\begin{array}{c}
\displaystyle
\sum_{n_1+n_2=n}
\frac{(-1)^{n_2}\Gamma(n+2+n_2+\mu)}{n_1!n_2!\Gamma(1+n_2+\mu)}\\
\\
\displaystyle
=\sum_{n_1+n_2=n}
\frac{(-1)^{n_2}}{n_1!n_2!}\Pi_{j=1}^{n+1}(j+n_2+\mu)\\
\\
\displaystyle
=\frac{1}{n!}\sum_{n_1+n_2=n}
\frac{(-1)^{n_2}n!}{n_1!n_2!}\left(\frac{d}{dx}\right)^{n+1}\{x^{n+1+n_2+\mu}\}\vert_{x=1}\\
\\
\displaystyle
=\frac{1}{n!}\left(\frac{d}{dx}\right)^{n+1}\{(1-x)^nx^{n+1+\mu}\}\vert_{x=1}\\
\\
\displaystyle
=\frac{n+1}{n!}\left(\frac{d}{dx}\right)^n\{(1-x)^n\}\cdot\frac{d}{dx}x^{n+1+\mu}\vert_{x=1}\\
\\
\displaystyle
=(-1)^n(n+1)(n+1+\mu).
\end{array}
$$

\noindent
$\Box$

\subsection{Proof of Lemma 4.1}

In \cite{Ik11} we have already known that
$$\begin{array}{c}
\displaystyle
v_{g_N}(y)
-v_{\tau}(y)\\
\\
\displaystyle
=-\sum_{m>N}\left\{\frac{(\tau-\sqrt{\tau^2+k^2})\overline\omega}{k}\right\}^mJ_m(kr)e^{im\theta}
-\sum_{m>N}\left\{\frac{(\tau+\sqrt{\tau^2+k^2})\omega}{k}\right\}^mJ_m(kr)e^{-im\theta},
\end{array}
\tag {A.6}
$$
where $y=(r\cos\,\theta,r\sin\,\theta)$.
Since
$$\begin{array}{c}
\displaystyle
\partial_{\tau}\left\{\frac{(\tau\mp\sqrt{\tau^2+k^2})\overline\omega}{k}\right\}^m
=m\left\{\frac{(\tau\mp\sqrt{\tau^2+k^2})\overline\omega}{k}\right\}^{m-1}
\frac{(\sqrt{\tau^2+k^2}\mp\tau)\overline\omega}{k\sqrt{\tau^2+k^2}}\\
\\
\displaystyle
=\mp\frac{m}{\sqrt{\tau^2+k^2}}\left\{\frac{(\tau-\sqrt{\tau^2+k^2})\overline\omega}{k}\right\}^{m},
\end{array}
$$
from (A.6) we have
$$\begin{array}{c}
\displaystyle
\sqrt{\tau^2+k^2}
(v_{\partial_{\tau}g_N}(y)
-\partial_{\tau}v_{\tau}(y))
\\
\\
\displaystyle
=\sum_{m>N}
m\left\{\frac{(\tau-\sqrt{\tau^2+k^2})\overline\omega}{k}\right\}^{m}
J_m(kr)e^{im\theta}
-\sum_{m>N}m\left\{\frac{(\tau+\sqrt{\tau^2+k^2})\omega}{k}\right\}^mJ_m(kr)e^{-im\theta}.
\end{array}
$$
Since
$$\displaystyle
\vert J_m(kr)\vert
\le\left(\frac{kr}{2}\right)^m\frac{1}{m!},
$$
it follows that
$$\displaystyle
\sqrt{\tau^2+k^2}
\vert
v_{\partial_{\tau}g_N}(y)
-\partial_{\tau}v_{\tau}(y)\vert
\le C(N+1)E(\tau;N+1),
\tag {A.7}
$$
where $C>0$ is independent of $N$ and $\tau$
and
$$\displaystyle
E(\tau;N)=\frac{1}{N!}
\left\{\frac{R(\tau+\sqrt{\tau^2+k^2})}{2}\right\}^N
e^{R(\tau+\sqrt{\tau^2+k^2})/2}.
$$
By virtue of the choice of $\tau(N)$, $\beta$ and the Stirling
formula (cf. \cite{O}), we have
$$\displaystyle
e^{R\tau(N)}E(\tau(N);N+1)=O(N^{-\infty}).
$$
This is the key point and see \cite{Ik11} for the detail of the
derivation.  This together with (A.7) yields the half of the
desired estimates.  The remaining estimate can be also given similarly.

\noindent
$\Box$

\end{document}